\documentclass[12pt]{article}

\usepackage{multicol}
\usepackage{amssymb}
\usepackage{epsfig}
\usepackage{amsmath}
\usepackage{graphics}

\usepackage{tikz}
\usepackage{mathrsfs,hyperref}
\usepackage{caption}

\usetikzlibrary{decorations.pathmorphing}
\usetikzlibrary{decorations.pathreplacing}
\usetikzlibrary{patterns}

\newenvironment{talign*}
 {\csname align*\endcsname}
 {\endalign}
\usepackage[utf8]{inputenc}
\usepackage[T1]{fontenc}

\usepackage{verbatim}
\usepackage{amsthm}
\usepackage{multicol, enumitem}

\newtheorem{thm}{Theorem}[section]
\newtheorem{cor}[thm]{Corollary}

\newtheorem{prop}[thm]{Proposition}

\newtheorem{ques}[thm]{Question}
\newtheorem{rem}[thm]{Remark}

\def\A{\mathcal{A}}

\def\F{\mathbb{F}}
\def\HH{\mathcal{H}}

\def\L{\mathcal{L}}

\def\P{\mathcal{P}}
\def\Q{\mathbb{Q}}

\def\a{\boldsymbol{a}}
\def\b{{\bf b}}
\def\c{{\bf c}}

\def\m{{\bf m}}

\def\k{{\bf k}}

\def\v{{\bf v}}
\def\w{{\bf w}}
\def\x{{\bf x}}
\def\y{{\bf y}}
\def\z{{\bf z}}
\def\0{{\bf 0}}
\def\1{{\bf 1}}

\begin{document}


\title{\Large {{\bf On Pappus Configurations in Hall Planes}}}
\bigskip

\bigskip

\author{{\bf Felix Lazebnik, Lorinda Leshock}\\
\small Department of Mathematical Sciences\\
\small University of Delaware, Newark, DE 19716, USA \\
\small {\tt fellaz@udel.edu}, {\tt
lorinda.leshock@gmail.com}}

\maketitle

\vspace{-0.2in}

\begin{abstract}
As the finite Hall planes are Non-Desarguesian, the Pappus Theorem does not hold in them. In this paper we state and prove some weaker versions of Pappus's Theorem in Hall planes.
\end{abstract}

\section{Introduction}\label{sectionIntroduction} Let $\P$ be a set whose elements we will call points, and $\L$ be a collection of subsets of $\P$ which we call lines. If a point is an element of a line, we say that it is {\it on} the line, or that the line {\it passes} through the point. In this situation, we may also say that the line is {\it on} the point, that the line contains the point, or that the point and line are incident. A set of points is called {\it collinear} if all points from the set are on the same line. A set of lines is called {\it concurrent} if all lines from the set are on the same point.
 We say that the pair $(\P, \L)$ is a {\it partial plane} or a {\it configuration}, if every two distinct points are on at most one line, and
every line contains at least two points. We say that configuration $(\P, \L)$ is isomorphic to a configuration $(\P', \L')$ if there exists a bijection $\P \to \P'$ such that the induced map $\L \to \L'$ is also a bijection. We say that a configuration $(\P, \L)$ is {\it embedded} in a configuration $(\P', \L')$ if there exists an injective map $\phi: \P \to \P'$ such that the image of every line $\ell \in \L$, defined as $\{\phi(P): P \in \ell\}$, is a subset of some line of $\L'$. When $(\P, \L)$ is embedded in $(\P', \L')$, we will also say that $(\P, \L)$ is {\it in} $(\P', \L')$.

We find the question whether or not a given configuration is embedded in a finite affine or a finite projective plane of great interest. Often the question is asked in the case when a partial plane is an affine or a projective plane itself. Some related results can be found in
Moorhouse and Williford
\cite{moorhouseWillifordPaper}, Lazebnik, Mellinger, and Vega \cite{lazebnikMellingerVegaPaper}, Metsch \cite{Metsch91}, Galiskan and Petrak \cite{caliskanPetrakPaper}, Petrak \cite{petrakPaper, petrakThesis}, Caliskan and Moorhouse \cite{caliskanMoorhousePaper}, Tait \cite{taitPaper}, and numerous references therein.

A celebrated result of Ostrom \cite{ostromPaper} establishes the existence of the Desargues configuration in every finite projective plane. The short proof in \cite{ostromPaper} is a beautiful pigeonhole argument that actually demonstrates the existence of many Desargues configurations in the following strong form:
\begin{quote}
	{\it In a finite projective plane, let $\ell_1$, $\ell_2$, and $\ell_3$ be three distinct lines through a point $P$. Let $R$ and $S$ be any two points not on $\ell_1$, $\ell_2$, or $\ell_3$. Consider the set of triangles with one vertex each on $\ell_1$, $\ell_2$, and $\ell_3$, one side going through $R$ and the other side going through $S$. At least one pair of triangles of this set satisfies Desargues's Theorem.}

\end{quote}
Finiteness of the plane is important: Hall's ``free plane'' construction in \cite{hallPaper} provides examples of infinite projective planes that do not contain the Desargues configuration.

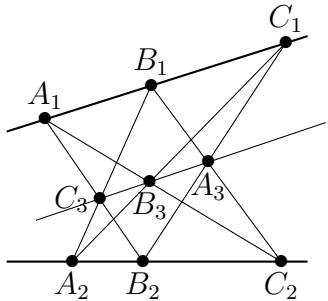
\begin{figure}[!h]
	\begin{center}
		\begin{tikzpicture}
		\draw[thick] (2,1) -- (6,1);
		\draw[thick] (2,2.73) -- (6,4);
		\draw[thin] (2.39,1.55) -- (6.24,2.85);
		\draw[thin] (2.87,1) -- (3.92,3.34) node {\textbullet};
		\draw[thin] (3.92,3.34) -- (5.65,1) node {\textbullet};
		\draw[thin] (2.51,2.89) -- (3.81,1) node {\textbullet};
		\draw[thin] (5.65,1) -- (2.51,2.89) node {\textbullet};
		\draw[thin] (3.81,1) -- (5.71,3.91) node {\textbullet};
		\draw[thin] (5.71,3.91) -- (2.87,1);
		\draw (4.68,2.32) node {\textbullet};
		\draw (3.9,2.06) node {\textbullet};
		\draw (3.24,1.83) node {\textbullet};
		\draw (2.87,1) node {\textbullet};
		\draw (2.51,2.89) node [above] {$A_1$};
		\draw (3.92,3.34) node [above] {$B_1$};
		\draw (5.71,3.91) node [above] {$C_1$};
		\draw (2.87,1) node [below] {$A_2$};
		\draw (3.81,1) node [below] {$B_2$};
		\draw (5.65,1) node [below] {$C_2$};
		\draw (3.24,1.83) node [left] {$C_3$};
		\draw (3.9,2.06) node [below] {$B_3$};
		\draw (4.68,2.32) node [below] {$A_3$};
		\end{tikzpicture}
	\end{center}
	\caption{Diagram of the Pappus configuration in a projective plane.}
\label{firstFigure}
\end{figure}

As far as the authors know, the existence of a Pappus configuration (see Figure \ref{firstFigure}\,) in every finite projective plane remains unknown. We believe that it does exist, and this paper grew from our failed attempt to prove this.

When a ternary ring that coordinatizes the projective plane ``contains'' a finite field, the Pappus configuration exists. Indeed, restricting coordinates of points $A_i,B_i,C_i$, $i=1,2$, to a finite field, leads to a Pappus configuration. The authors are not aware of any example of a finite projective plane for which every coordinatization with a ternary ring ``does not contain'' a finite field. Even if such a plane exits, it still may contain Pappus configurations.

Moreover, in all planes we considered, the number of Pappus configurations was actually large. This was not too surprising, since the collineation group of the plane acts on the set of its Pappus configurations. Therefore, when
the existence of Pappus configurations in a plane was clear, we tried to find a statement as close as possible to Pappus's Theorem in the plane.
Here is the motivation for the main result of this paper. It would imply the existence of Pappus configurations in finite projective planes in a strong way.

\begin{ques}[$3+2$ Question]\label{ques32} Is it true that in a finite projective plane the following holds: For every pair of lines $\ell_1, \ell_2$ and every three points on $\ell_1$, and every two points on $\ell_2$, one more point on $\ell_2$ can be found so that the six points define a Pappus configuration?
\end{ques}

A question with a weaker condition is the following.

\begin{ques}[$3+1$ Question]\label{ques31} Is it true that in a finite projective plane the following holds: For every pair of lines $\ell_1, \ell_2$ and every three points on $\ell_1$, and every one point on $\ell_2$, two more points on $\ell_2$ can be found so that the six points define a Pappus configuration?
\end{ques}

The $3+1$ Question was answered affirmatively for all planes of order less than $25$, and for some of order $25$ using the data available on Moorhouse's database of projective planes \cite{moorhouseWebsite}\,, or using built-in Magma commands and our own code. In particular, it is answered affirmatively for Hall planes of orders $9, 16, 25$; for the Hughes plane of order $25$; for Czerwinski \& Oakden planes of order $25$: $a1, a6, b3, b6$; and for Rao planes of order $25: a5, a7$. As of yet, we have been unable to answer the $3+1$ Question completely even for Hall planes of orders greater than $25$.

It is easy to argue that the $3+1$ Question is affirmed in all affine planes of order $n$ if and only if it is affirmed in all projective planes of order $n$.

The $3+2$ Question was answered negatively for all planes of order $25$ and for some planes of order $49$: Hall plane, Hughes plane, and Dickson Near-field plane. Moreover, we have not found finite nonclassical planes (i.e., the ones that cannot be coordinatized by a field) where the $3+2$ Question is answered affirmatively. Clearly, it takes longer (for the same plane order) to verify numerically that the $3+1$ Question has the affirmative answer than to obtain a negative answer to the $3+2$ Question. That is why the latter could be tested for larger planes.

If one strengthens the condition of the $3+2$ Question requiring that three arbitrary points on line $\ell_2$ are chosen, a similar ``$3+3$ Question'' is equivalent to Pappus's Theorem, and hence, the answer is affirmative only in classical planes. One can weaken the condition of the $3+1$ Question further by not requiring that any arbitrary points on line $\ell_2$ are chosen and call it a ``$3+0$ Question''.

In order to state our results, we need additional definitions and some preliminary results. They are collected in Section \ref{sectionDefinitions}\,.
In Section \ref{sectionResults}\,, we describe our results for Hall planes, and their proofs appear in Section \ref{sectionProofs}\,. In Section \ref{sectionConcluding}\,,
we make concluding remarks and mention several open problems.

\bigskip

\section{Definitions, notations, and preliminary results.} \label{sectionDefinitions} Our exposition is based on Hall \cite{hallPaper, hallCombiBook}, and Leshock \cite{leshock}. In \cite{leshock}, many proofs that were not explicit in \cite{hallPaper, hallCombiBook} were checked analytically.

\subsection{Hall system}\label{subsectionHallSystem}
Let $\F=F_q$ represent the finite field of prime power order $q$ and $\HH=\{ (a_1,a_2) \, : \, a_1, a_2 \in \F \}$. We call $\F$ the {\it basefield} of $\HH$. Clearly, $|\HH|=q^2$. It is convenient to have multiple notations for elements of $\HH$. Let the bold letter $\a$ denote the ordered pair, $(a_1,a_2) \in \HH$. We will identify an element $\b = (b_1,0)$ of $\HH$ with the element $b_1$ of $\F$ and write, for simplicity, $\b \in \F$, and we will write $\b=(b_1,b_2) \not \in \F$ when $b_2 \not=0$.

The {\it Hall system} is a two dimensional (right and left) vector space $\HH$ over $\F$ equipped with a certain multiplication of vectors. The addition in the Hall system $\{\HH,+,\boldsymbol{\cdot} \}$ is the usual addition in $\F^2$. To define multiplication in $\{\HH, +, \boldsymbol{\cdot} \}$, we use the operations from the basefield and a quadratic polynomial $f(x)=x^2-rx-s$ with $r,s \in \F$ which is irreducible over $\F$, and we refer to it as the {\it defining polynomial} of the system. More specifically, for $\a,\b \in \HH$,
\begin{itemize}
\item[$A.$] $\a + \b=(a_1, a_2) + (b_1, b_2) = (a_1+b_1, a_2+b_2)$
\item[$M1.$] $\a \b = (a_1, a_2) \boldsymbol{\cdot} (b_1, b_2) =(a_1 b_1, a_2 b_1)$, if $\b \in \F$
\item[$M2.$] $\a \b = (a_1, a_2) \boldsymbol{\cdot} (b_1, b_2) = (a_1 b_1-a_2 b_2^{-1}f(b_1), a_1 b_2-a_2 b_1+a_2 r)$, if $\b \not \in \F$
\end{itemize}

For Hall systems over basefields with $q>3$, multiplication is neither commutative nor associative; it is right distributive but not left distributive over addition. Clearly, $\HH$ is a group under addition with identity element $\0 = (0,0)$. It's also easy to see that the identity element for multiplication in $\HH$ (both left and right) is $\1=(1,0)$.

\subsection{Hall affine plane}\label{subsectionHallAffine}
Having a Hall system $\HH$, we can construct a {\it Hall affine plane} $\A_\HH$ as follows. The point set of $\A_\HH$ is $\P = \left \{ ( \x,\y ) \, : \, \x, \y \in \HH \right \}$. For arbitrary fixed $\m,\k \in \HH$, the sets $\{(\x,\x\m+\k) \, : \, \x \in \HH \}$, represent ``non-vertical'' lines with equations of the form $\y=\x\m+\k$. For arbitrary fixed $\k \in \HH$, the sets $\{(\x,\k) \, : \, \x \in \HH \}$, represent ``horizontal'' lines with equations of the form $\y=\k$. For arbitrary fixed $\c \in \HH$, the sets $\{(\c,\y) \, : \, \y \in \HH \}$, represent ``vertical'' lines with equations of the form $\x=\c$. All these $q^4 +q^2$ lines form the set $\L$ of lines $\A_\HH$. It is easy to verify that $\A_\HH = (\P,\L)$ is an affine plane of order $q^2$. The non-vertical lines $\y = \x \m + \k$ will be referred to as {\it type 1} lines or {\it type 2} lines depending on the ``slope'' $\m$; type $1$ lines have $\m \in \F$, and type $2$ lines have $\m \not \in \F$.

\subsection{Collineation of Hall affine plane}\label{subsectionCollineation}
Our exposition here is based on Hughes \cite{hughesCollPaper2}. Understanding the action of the collineation group of $\A_\HH$ on points, lines, and pairs of lines of $\A_\HH$ was crucial for our work on the $3+1$ Question, and some facts presented in this section are original. Those proofs which we omit (many are straightforward) can be found in \cite{hughesCollPaper2} and \cite{leshock}. All notions and facts related to groups and group actions, used in this section, can be found in \cite{hallGroupBook}.

Let $\A = (\P, \L)$ be an affine plane. A bijection $\phi: \P \to \P$ which preserves collinearity of points in $\A$ is called a {\it collineation}. A collineation $\phi$ acts on the set of lines $\L$ in an obvious way: for $\ell \in \L$, $\phi(\ell) = \{\phi( P) \, : \, P \in \ell \}$. A point $P$ (line $\ell$) is {\it fixed} under $\phi$ if $\phi(P) = P$ ($\phi(\ell) = \ell$). We say that a line $\ell$ is {\it fixed pointwise} when every point on $\ell$ is fixed. The parallel class of lines is {\it fixed} under a collineation if the collineation permutes all lines in the class.

A {\it translation} of an affine plane $\A$ is a collineation of $\A$ such that the parallel classes are fixed and there is a parallel class such that every line from the class is fixed. It is easy to see that for fixed $\a,\b\in \HH$, the map $\tau_{\a,\b}: \P \to \P$ given by $\tau_{\a,\b}\left((\x,\y)\right) = (\x +\a,\y+\b)$ is a translation.

Let $TR = \{ \tau_{\a,\b} \, : \, \a, \b \in \HH \}$ be the set of all translations $\tau_{\a,\b}$ of $\A_\HH = (\P, \L )$. For arbitrary fixed $\c,\m, \k \in \HH$, following \cite{hughesCollPaper2} and Kallaher \cite{kallaherChapter}, we will denote the vertical line $\{(\c,\y) \, : \, \y \in \HH \}$ as $[\c]$, and the non-vertical line $\{(\x,\x\m+\k) \, : \, \x \in \HH \}$ as $[\m,\k]$.
Then it's easy to check that $\tau=\tau_{\a,\b}$ acts on $\L$ in the following way:
\[
\tau\left([\c]\right) = [\c+\a];\;\;\;\;\;\;\tau\left([\m, \k]\right) = [\m, \k - \a \m + \b].
\]
We summarize the properties of $TR$ below.

\begin{prop} \label{2.1.2} $TR$ is a subgroup of the collineation group of $\A_\HH$ and it is isomorphic to the additive group of the vector space $\HH^2$ over $\F$. It is sharply transitive on $\P$,
preserves the classes of parallel lines and acts transitively on the set of lines in every parallel class.
\end{prop}

If the group of translations of an affine plane $\A$ is transitive on the points of $\A$, the plane $\A$ is called a {\it translation plane}, see \cite{ostromTransConfigs}. Thus, the Hall affine plane $\A_\HH$ is a translation plane. It is known that translation planes are exactly the planes which can be coordinatized by quasifields (which we don't need to define), of which the Hall system is a particular example.

In $1959$, Hughes gave an analytic description of the entire collineation group of Hall planes, see \cite{hughesCollPaper2} and his related work \cite{hughesCollPaper3}.

\begin{rem}\label{hughesRemark} {\rm Hughes used different conventions than Hall did in his construction of Hall planes, e.g., multiplication is left not right distributive over addition and the equation of a non-vertical line is of the form $\m \x + \y = \k$ rather than $\y = \x \m+ \k$. We express Hughes' collineation subgroups using Hall's conventions.}
\end{rem}

Hughes presents the collineation group of Hall planes as generated by six subgroups. Three of those subgroups are relevant to this paper. One of them, $TR$ has been already presented, and the other two, ATP and LNR, are presented below.

Let $S$ be a $2 \times 2$ nonsingular matrix over $\F$. For $\x \in \HH$, let $S$ act on $\HH$ as
$\x S = (x_1, x_2) S.$
Define $\a = \1 S^{-1}$, where $\1 = (1,0)$ is the multiplicative identity of $\HH$. Hence, $\a S =\1$.

Let $\sigma_S$ be a map $\sigma_S: \P \to \P$, $\sigma_S\left((\x,\y)\right) = (\x S,\y S)$. Such a map $\sigma_S$ is called an
{\it autotopism} of $\A_\HH$, and it is easy to check that it is a collineation of $\A_\HH$.
Let $ATP = \{\sigma_S \, : \, S \in GL(2,\F) \}$ be the set of all autotopisms of $\A_\HH = (\P, \L)$. Then it's easy to check that $\sigma=\sigma_S$ acts on $\L$ in the following way,
\[\sigma\left([\c]\right) = [\c S];\;\;\;\;\;
\sigma \left([\m, \k] \right) = [(\a \m )S, \k S].
\]

We summarize the properties of autotopisms as follows.

\begin{prop}\label{2.1.3} $ATP$ is a subgroup of the collineation group of $\A_\HH$,
and it is isomorphic to $GL(2, \F)$.
 It fixes the parallel classes of type $1$ lines, acts transitively on the parallel classes of type $2$ lines, and it has two orbits on the set of vertical lines $\left( \text{one orbit consists of the line } [\0] \right)$.
\end{prop}

Finally, we present the third group of collineations of $\A_\HH$ that we use.
Let $a,b \in \F$ not both zero, and define the matrix $L=\begin{pmatrix} -a r + b &a s \\ a & b \end{pmatrix}$, where $r,s$ are the coefficients of the defining polynomial $f(x)=x^2-r x-s$ of $\HH$. Clearly, $L \in GL(2, \F)$.

\begin{rem} {\rm Our matrix $L$ is the transpose and has a different sign in position $(1,1)$ than the corresponding matrix in \cite{hughesCollPaper3}. This is due to the adjustment of Hughes' conventions to Hall's conventions, see Remark \ref{hughesRemark}\,.}
\end{rem}

For $(\x,\y) \in \HH^2$, let $L$ act on $\HH^2$ as $(\x,\y) L=((-ar+b)\x+a \y, as \x+b \y).$
Let $\lambda=\lambda_L=\lambda_{a,b}$ be a map $\lambda: \P \to \P$ given by $\lambda \left((\x,\y) \right) = (\x,\y)L.$
Then it's easy to check that $\lambda$ acts on $\L$ in the following way:
\begin{talign*}
\lambda \left([\c] \right) &= \begin{cases}
[b \c ], & \text{ if } a=0 \\
\left[ \left(\frac{b}{a},0 \right), - \frac{b^2-abr-a^2s}{a}\c \right], &\text{ if } a \not=0
\end{cases} \\
\text{If $\m \in \F$,} \\
\lambda \left( [\m, \k] \right) &= \begin{cases}
[\m, b \k], &\text{ if } a = 0 \\
[a \k ], &\text{ if } a \not=0 \text{ and } m_1 =r-\frac{b}{a} \\
\left[ \left(\frac{a s+m_1 b}{a m_1-ar+b},0\right),\frac{b^2-abr-a^2s}{a m_1-ar+b}\k \right], &\text{ if } a \not= 0 \text{ and } m_1 \not=r-\frac{b}{a}
\end{cases} \\
\text{If $\m \not \in \F$,} \\
\lambda \left([\m, \k] \right) &=[\m, -a \k \m +b \k ]
\end{talign*}
Such a map $\lambda$ is called a {\it linear} map of $\A_\HH$. Let $LNR = \left\{\lambda_{a,b} \, : \, a,b \in \F, (a,b) \not= (0,0) \right\}$ be the set of all linear maps of $\A_\HH = (\P, \L )$.

We will denote the set of lines that are of type $1$ or vertical as {\it BF} (slope is from \underline{b}ase\underline{f}ield or $\infty$), and the set of lines that are of type $2$ as {\it NBF} (slope is \underline{n}ot from the \underline{b}ase\underline{f}ield). There are fewer $BF$ lines than $NBF$ lines: $q^3+q^2$ of the former and $q^4-q^3$ of the latter.

\begin{prop} \label{2.1.4}
 $LNR$ is a subgroup of the collineation group of $\A_\HH$, and it is isomorphic to the multiplicative group of the quadratic extension field $\F[\alpha]=F_q[\alpha]$ where $\alpha^2=-r \alpha+s$, and so to $F_{q^2}^\times$. $LNR$ fixes the parallel classes of lines in $NBF$ and is transitive on the set of all parallel classes of lines in $BF$.
\end{prop}

We conclude this subsection with several statements describing the action of the groups
$TR$, $ATP$ and $LNR$ on lines of $\A_\HH$.

\begin{prop}\label{tratplnr} (i) The group generated by $TR$ and $ATP$ acts transitively on all type $2$ lines and on all vertical lines.

(ii) The group generated by $TR$ and $LNR$ acts transitively on all $BF$ lines.
\end{prop}

\begin{proof} (i) By Proposition \ref{2.1.3}\,, $ATP$ acts transitively on the parallel classes of type $2$ lines and has two orbits on the set of vertical lines, one of which is $[\0]$. By Proposition \ref{2.1.2}\,, $TR$ acts transitively on the lines in every parallel class.

(ii) By Proposition \ref{2.1.4}\,, $LNR$ acts transitively on the parallel classes of lines in $BF$. By Proposition \ref{2.1.2}\,, $TR$ acts transitively on the lines in every parallel class.
\end{proof}

Based on the above propositions, it is clear that there are at most two orbits of lines in $\A_\HH$. The sets $BF$ and $NBF$ partition $\L$ and are precisely these orbits. This was implicitly shown in \cite{hughesCollPaper2}.

The next proposition describes the action of certain collineations in the stabilizer of the origin on lines through the origin.

\begin{prop}\label{oneLine} \cite{leshock} For any line through the origin, there exists a group of collineations that fixes the line, fixes the origin, and acts transitively on the other points of the line.
\end{prop}

\begin{proof} We consider the following three cases: (i) $\x = \0$, (ii) $\y = \x \m \text{ with } \m \in \F$, and (iii) $\y = \x \m \text{ with } \m \not\in \F.$

(i) Take $\ell: \x = \0$. The line $\ell$ is on the origin. Every point of $\ell$ is of the form $(\0,\y)$ for some $\y \in \HH$. If $\y\not= (0,0)$, then for every nonzero $\z \in \HH$, there exists $S \in GL(2,\F)$ such that $ \y S=\z$. Hence, $ATP$ acts transitively on points of the $\ell: \x=\0$ distinct from the origin.

(ii) Take $\ell: \y = \x \m, \text{ with } \m \in \F$. The line $\ell$ is on the origin. Every point of the line $\ell$ is of the form $(\x, \x \m)$ for $\x \in \HH$. If $(\x, \x\m)$ and $(\z,\z\m)$ are two nonzero points of $\ell$, then the same $S$ as we used in the proof of part (i), will fix $\ell$ and map $(\x, \x\m)$ to $(\z,\z\m)$. Hence, $ATP$ acts transitively on points of the $\ell$ distinct from the origin.

(iii) Take $\ell: \y = \x \m, \text{ with } \m \not\in \F$. The line $\ell$ is on the origin. By Proposition \ref{2.1.4}\,, the action of $LNR$ fixes the parallel class of line $\ell$. Clearly, $LNR$ fixes the origin. Hence, it fixes the line $\ell$.

Let $(\v, \v \m)$ and $(\w, \w \m)$ be two distinct nonzero points of $\ell$. Let us show that there exists $a,b \in \F$ with $(a,b) \not=(0,0)$ such that $\lambda= \lambda_{a,b} \in LNR$ maps $(\v,\v \m)$ to $(\w,\w \m)$. This is equivalent to solving the equation $(-a r+b)\v+a \v \m = \w$ for $a,b$.

Since $\HH$ is a two dimensional vector space over $\F$, and $\v$ and $\v \m$ are linearly independent over $\F$ because $\m \not \in \F$, there exist $c_1, c_2 \in \F$ such that $c_1 \v+ c_2 \v \m = \w$. Setting $-ar+b=c_1$ and $a= c_2$, we get $b = c_1+c_2 r$. Clearly, $(a,b) \not =(0,0)$ as otherwise, $\w = \0$. Hence, the transitivity statement is proven.
\end{proof}

\begin{cor} \cite{leshock} \label{fixPtCor} (i) Any pair of intersecting lines can be mapped by a translation of $\A_{\HH}$ to a pair of lines $\{ \ell_1^\prime, \ell_2^\prime\}$ meeting at the origin. Furthermore, if $\ell_1^\prime$ is in $BF$, then there exists a collineation which maps it to $\ell_1: \x=\0$. Otherwise, $\ell_1^\prime$ is in $NBF$, and there exists a collineation which maps it to $\ell_1: \y = \x (0,1)$. Let $\ell_2$ be the image of the second line under either of these two maps. There exists a subgroup of the collineation group of the plane which fixes the origin, fixes lines $\ell_1$ and $\ell_2$, and acts transitively on the points distinct from the origin of $\ell_1$, or acts transitively on the points distinct from the origin of $\ell_2$.

(ii) Any pair of parallel lines can be mapped by a translation of $\A_{\HH}$ to a pair of lines $\{ \ell_1^\prime, \ell_2^\prime\}$ with $\ell_1^\prime$ on the origin. Furthermore, if $\ell_1^\prime$ is in $BF$, then there exists a collineation which maps it to $\ell_1: \y=\0$. Otherwise, $\ell_1^\prime$ is in $NBF$, and there exists a collineation which maps it to $\ell_1: \y = \x (0,1)$. Let $\ell_2$ be the image of the second line under either of these two maps. There exists a subgroup of the collineation group of the plane which fixes the origin, fixes line $\ell_1$, fixes the parallel class of line $\ell_2$, and acts transitively on the points distinct from the origin of $\ell_1$, or acts transitively on the points of $\ell_2$.
\end{cor}

\begin{proof} (i) The first statement follows directly from Proposition \ref{2.1.2}\,. Since $\ell_1^\prime$ and $\ell_2^\prime$ intersect at the origin, any collineation that fixes the origin and the parallel class of each line, fixes both lines.
\begin{itemize}
\item[Case $1$]: If $\ell_1^\prime, \ell_2^\prime \in BF$, since $LNR$ fixes the origin, $\ell_1^\prime$ can be mapped to $\ell_1: \x=\0$ by Proposition \ref{2.1.4}\,. Then the image of $\ell_2^\prime$ is $\ell_2 : \y = \x \m$ for some $\m \in \F$ under this mapping. Following part (i) of the proof of Proposition \ref{oneLine}\,, we conclude that there exists a map in $ATP$ that acts transitively on the points of $\ell_1$ different from the origin. By Proposition \ref{2.1.3}\,, and the fact that elements of $ATP$ fix the origin and the pair of lines intersect at the origin, it maps $\ell_2$ to itself. Alternatively, following part (ii) of the proof of Proposition \ref{oneLine}\,, we conclude that there exists a map in $ATP$ that acts transitively on the points of $\ell_2$ different from the origin. By Proposition \ref{2.1.3}\,, and the fact that elements of $ATP$ fix the origin and the pair of lines intersect at the origin, it maps $\ell_1$ to itself.

\item[Case $2$]: If $\ell_1^\prime, \ell_2^\prime \in NBF$, we first use a collineation from $ATP$ to map $\ell_1^\prime$ to $\ell_1: \y=\x(0,1)$. Such exists by Proposition \ref{2.1.3}\,. Let $\ell_2^\prime$ be mapped to a line $\ell_2$ by this map. Then $\ell_1$ and $\ell_2$ intersect at the origin and are in $NBF$. By part (iii) of the proof of Proposition \ref{oneLine}\,, $LNR$ acts transitively on points of $\ell_1$ distinct from the origin, fixes the origin and the slope of any line in $NBF$. Hence, it fixes $\ell_2$. Alternatively, we can show that there exists a collineation in $LNR$ that acts transitively on points of $\ell_2$ distinct from the origin, fixes the origin and the slope of any line in $NBF$. Hence, it fixes $\ell_1$.

\item[Case $3$]: If $\ell_1^\prime \in BF$ and $\ell_2^\prime \in NBF$, then we can use a collineation of $LNR$ to map $\ell_1^\prime$ to $\ell_1:\x=\0$, by Proposition \ref{2.1.4}\,. Next, by Proposition \ref{2.1.3}\,, we can use a collineation of $ATP$ to map $\ell_2^\prime$ to $\ell_2: \y = \x(0,1)$. Note that this collineation necessarily fixes $\ell_1$.

Now, we show that for every two nonzero points on $\ell_1$, $(\0,\y)$ and $(\0,\z)$, where $\y = (y_1, y_2)$ and $\z = (z_1,z_2)$, there exists a nonsingular $2 \times 2$ matrix $S$ over $\F$ such that $\y S = \z$ and $\left(\left(\1 S^{-1}\right) (0,1)\right) S = (0,1)$. This will ensure that there is a subgroup of $ATP$ that fixes both lines, fixes the origin, and acts transitively on the points of $\ell_1$ distinct from the origin. Searching for $S$ with four undetermined entries, we obtain that 
$$S=\begin{pmatrix} \frac{y_1 z_1 + r y_2 z_1 - s y_2 z_2}{y_1^2 + r y_1 y_2 - s y_2^2} & \frac{-y_2 z_1 + y_1 z_2}{y_1^2 + r y_1 y_2 - s y_2^2} \\ \frac{s (-y_2 z_1 + y_1 z_2)}{y_1^2 + r y_1 y_2 - s y_2^2} & \frac{y_1 z_1 + r y_1 z_2 - s y_2 z_2}{y_1^2 + r y_1 y_2 - s y_2^2}\end{pmatrix}.$$
 The denominators of fractions of entries of $S$ are nonzero, due to points $\y, \z$ being nonzero, and the irreducibility of the defining polynomial $f(x)=x^2 - rx -s$ of $\HH$ over $\F$.

Similarly, we can show that there exists a subgroup of $ATP$ which fixes both lines, fixes the origin, and acts transitively on the points of $\ell_2$ distinct from the origin.
\end{itemize}

(ii) As lines are parallel, then both are either in $BF$ or in $NBF$. The proof follows the same ideas as the one of part (i) and is shorter. Because of this, we omit it.
\end{proof}

Now, we consider the action of the collineation group of $\A_\HH$ on the set of all ordered pairs of distinct lines which we refer to as just {\it pairs of lines}. As the collineation group acts on Pappus configurations, our proofs of the main results can be restricted to Pappus configurations on special pairs of lines. Since we have two orbits of lines, $BF$ and $NBF$, we consider three cases corresponding to choosing two lines from $BF$, or from $NBF$, or one line from each orbit. Propositions \ref{bfbfCases}\,, \ref{nbfnbfCases}\,, and \ref{nbfbfCases}\, below follow easily from Corollary \ref{fixPtCor}\, and will be instrumental for the proofs of the main results in Section \ref{sectionProofs}\,.

\begin{prop}[$BF$/$BF$] \label{bfbfCases}
(i) Any pair of intersecting lines from the $BF$ orbit can be mapped by a collineation of $\A_{\HH}$ to a pair of lines $(\ell_1, \ell_2)$ where $\ell_1: \y = \x \boldsymbol{\mu}$ for some $\boldsymbol{\mu} \in \F$ and $\ell_2: \x = \0$. Moreover, for any point on $\ell_1: A_1 ((\alpha,\beta), (\mu_1 \alpha, \mu_1 \beta))$, such a map can be found so that $A_1$ is mapped to $((0,1),(0,\mu_1))$.

(ii) Any pair of parallel lines from the $BF$ orbit can be mapped by a collineation of $\A_{\HH}$ to a pair of horizontal lines $(\ell_1, \ell_2)$ where $\ell_1: \y=\boldsymbol{\kappa}$ and $\ell_2: \y=\0$ for some $\boldsymbol{\kappa} \in \HH$. Moreover, for any point on $\ell_1: A_1 ((\alpha,\beta), (\kappa_1,\kappa_2))$, such a map can be found so that $A_1$ is mapped to $((0,1),(\kappa_1,\kappa_2))$.
\end{prop}

\begin{prop}[$NBF$/$NBF$] \label{nbfnbfCases}
(i) Any pair of intersecting lines from the $NBF$ orbit can be mapped by a collineation of $\A_{\HH}$ to a pair of lines $(\ell_1,\ell_2)$ where $\ell_1: \y = \x (\mu, \psi)$ for some $\mu, \psi \in \F$ and $\ell_2: \y = \x (0,1)$. Moreover, for any point on $\ell_1: A_1 \left( \left(\alpha,\beta \right), \left(\alpha \mu - \frac{\beta(-s-r \mu+\mu^2)}{\psi},r \beta-\beta\mu +\alpha \psi \right) \right)$, such a map can be found so that $A_1$ is mapped to $\left( \left(0,1 \right), \left(\frac{-s-r \mu+\mu^2}{\psi},r-\mu \right) \right)$.

(ii) Any pair of parallel lines from the $NBF$ orbit can be mapped by a collineation of $\A_{\HH}$ to a pair of lines $(\ell_1,\ell_2)$ where $\ell_1: \y=\x(0,1)+\boldsymbol{\kappa}$ and $\ell_2: \y=\x(0,1)$ for some $\boldsymbol{\kappa} \in \HH$. Moreover, for any point on $\ell_1: A_1 ((\alpha,\beta), (\kappa_1+s \beta, \kappa_2+\alpha+r \beta))$, such a map can be found so that $A_1$ is mapped to $((0,1),(\kappa_1+s,\kappa_2+r))$.
\end{prop}

Due to the asymmetry in the versions of Pappus's theorem that we will establish, the role of lines $\ell_1$, $\ell_2$ in the ordered pair $(\ell_1,\ell_2)$ is not symmetric either. Therefore, we need to consider two cases: line $\ell_1$ from $NBF$ and line $\ell_2$ from $BF$, and line $\ell_1$ from $BF$ and line $\ell_2$ from $NBF$. Also, note that if one line in a pair comes from $BF$ and another from $NBF$, they cannot be parallel.

\begin{prop}[$NBF$/$BF$ \& $BF$/$NBF$] \label{nbfbfCases}
(i) Any pair of intersecting lines with the first from $NBF$ and the second from $BF$ can be mapped by a collineation of $\A_{\HH}$ to $\ell_1: \y = \x (0,1), \ell_2: \x = \0$. Moreover, for any point on $\ell_1: A_1 ((\alpha,\beta), (s \beta, \alpha + r \beta))$, such a map can be found so that $A_1$ is mapped to $((0,1),(s,r))$.

(ii) Any pair of intersecting lines with the first from $BF$ and the second from $NBF$ can be mapped by a collineation of $\A_{\HH}$ to $\ell_1: \x = \0, \ell_2: \y = \x (0,1)$. Moreover, for any point on $\ell_1: A_1 ((0,0), (\alpha,\beta))$, such a map can be found so that $A_1$ is mapped to $((0,0),(0,1))$.
\end{prop}

\bigskip

\section{Results}\label{sectionResults} In Section \ref{sectionIntroduction}\, we mentioned the $3+0$ Question which motivates our results. We also presented the $3+1$ Question together with some comments. Here we wish to mention several partial results concerning the $3+1$ Question.

Recall that in the Hall affine plane of order $q^2$, the point set consists of ordered pairs of elements of a Hall system, and that each element of the Hall system may be represented as an ordered pair of elements of the basefield $F_q$, where $q$ is a prime power. We call each basefield element of a point a {\it component}. If a line is fixed, and we want to choose a point on it, it is sufficient to choose either the $\x$-coordinate or the $\y$-coordinate (in the case of a vertical line) of the point. It is clear that for fixed lines $\ell_1$ and $\ell_2$, the number of ways of choosing an ordered triple of points on $\ell_1$ and one point on $\ell_2$ is $\sim \left(q^2 \right)^3 q^2=q^8$, $q \to \infty$. When we can answer the $3+1$ Question affirmatively for all but $O(q^7)$ choices of the four points, we say the $3+1$ Question is affirmed {\it asymptotically}.

We succeeded affirming the $3+1$ Question asymptotically in some cases, but not in all.
\begin{thm}\label{theo31} The $3+1$ Question is affirmed asymptotically in $\A_\HH$ when two lines $\ell_1, \ell_2$ are both from $BF$ or are both from $NBF$.
\end{thm}

We succeeded affirming the $3+0$ Question in more cases.
\begin{thm}\label{theo30} The $3+0$ Question is affirmed in $\A_\HH$ when two lines $\ell_1, \ell_2$ are both from $BF$ or are both from $NBF$ and when line $\ell_1$ is from $NBF$ and line $\ell_2$ is from $BF$.
\end{thm}

How can one compare an asymptotic $3+1$ result to a complete $3+0$ result on a given pair of lines? We don't see an obvious way to compare them, and we think that is a matter of taste of a reader.

Our most general result is an understatement of what we believe to be true. Here is our main result.

\begin{thm}[$2+0$ Theorem]\label{theo20} In a Hall plane the following holds: For every pair of lines $\ell_1, \ell_2$ and every two points on $\ell_1$, a third point on $\ell_1$ and three points on $\ell_2$ can be found so that the six points define a Pappus configuration.
\end{thm}

\bigskip

\section{Proofs of Theorems} \label{sectionProofs}

\begin{figure}[h!]
	\begin{center}
		\begin{tikzpicture}[scale=0.75]
\draw[black, thick] (14.00,2.14) -- (1.03,10.45) node[above left] {$\ell_2: [\boldsymbol{\psi},\boldsymbol{\rho}]$}; 
\draw[black, thick] (13.97,2.99) -- (0.58,3.02) node[above left] {$\ell_1: [\boldsymbol{\mu},\boldsymbol{\kappa}]$}; 
\draw[black, thick] (9.34,3.78) --(2.59,5.61) node[above left] {$\ell_3$}; 
\draw[black, thick] (1.57,10.11) -- (5.80,3.01); 
\draw[black, thick] (1.57,10.11) -- (7.65,3); 
\draw[black, thick] (7.09,6.57) -- (1.06,3.02); 
\draw[black, thick] (7.09,6.57) -- (7.65,3); 
\draw[black, thick] (8.91,5.40) -- (1.06,3.02); 
\draw[black, thick] (8.91,5.40) -- (5.80,3.01); 
\draw (1.57,10.11) node {\textbullet}; 
\draw (7.09,6.57) node {\textbullet}; 
\draw (8.91,5.40) node {\textbullet}; 
\draw (1.06,3.02) node {\textbullet}; 
\draw (5.80,3.01) node {\textbullet}; 
\draw (7.65,3.00) node {\textbullet}; 
\draw (12.67,2.99) node {\textbullet}; 
\draw (4.57,5.08) node {\textbullet}; 
\draw (6.28,4.60) node {\textbullet}; 
\draw (7.45,4.28) node {\textbullet}; 
\draw (1.57,10.11) node[above right] {$A_2 ((g,h),(g,h)\boldsymbol{\psi}+\boldsymbol{\rho})$}; 
\draw (7.09,6.57) node[above right] {$B_2 ((t,v),(\dots))$}; 
\draw (8.91,5.40) node[above right] {$C_2 ((w,z),(\dots))$}; 
\draw (1.06,3.02) node[below] {$A_1 ((\alpha,\beta),(\alpha,\beta)\boldsymbol{\mu}+\boldsymbol{\kappa})$}; 
\draw (5.80,3.01) node[below] {$B_1 ((\gamma, \delta),(\dots))$}; 
\draw (7.65,3.00) node[below right] {$C_1 ((e, j),(\dots))$}; 
\draw (12.67,2.99) node[above right] {$O$}; 
\draw (4.57,5.08) node[below] {$C_3$}; 
\draw (6.28,4.60) node[below] {$B_3$}; 
\draw (7.45,4.28) node[below] {$A_3$}; 
		\end{tikzpicture}
	\end{center}
	\caption{Diagram of Pappus configuration with labels for the $2+0$~Question, no line is vertical.}
\label{secondFigure}
\end{figure}
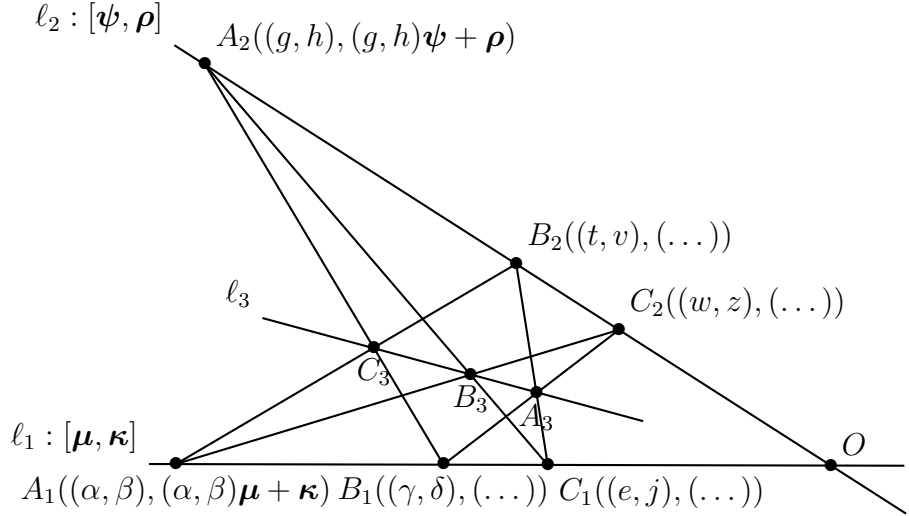

Let us describe the ideas used in the proofs of our results. As we study Hall planes analytically, our technique is ``just'' analytic geometry over Hall systems.
The proofs of existence of Pappus configurations follow from showing that certain systems of equations have solutions. In general, finding all solutions is not feasible. Therefore, we try to find at least one particular solution that will correspond to finding a Pappus configuration. Let us call the components of given points and the coefficients in the equations of lines $\ell_1$ and $\ell_2$ {\it parameters}, and the components of points whose existence we wish to establish {\it unknowns}. In all those cases, our approach is to specialize some unknowns and determine others. In all cases, we try to minimize the number of parameters by using the properties of the collineation group of $\A_\HH$. The details are stated in Propositions \ref{bfbfCases}\,, \ref{nbfnbfCases}\,, and \ref{nbfbfCases}\,. For example, to prove the $2+0$ Question is affirmative, we would assume that points $A_1$ and $B_1$ are given and points $C_1, A_2, B_2$, and $C_2$ are to be found. See Figure \ref{secondFigure}\,. To make it easier to distinguish between parameters and unknowns, we denote parameters by Greek letters and unknowns by Latin letters, except for $r$ and $s$, which will always denote the coefficients of the defining polynomial $f(x) = x^2-r x - s$ of $\HH$.

In Figure \ref{secondFigure}\,, the components of the given points $A_1$ and $B_1$ on line $\ell_1$ correspond to the sequence of parameters $(\alpha, \beta, \gamma, \delta)$. The elements of the Hall system representing the coefficients of the lines are $(\boldsymbol{\kappa}, \boldsymbol{\mu}, \boldsymbol{\rho}, \boldsymbol{\psi})$. The components of point $C_1$ on line $\ell_1$ and points $A_2, B_2$, and $C_2$ on line $\ell_2$, the existence of which we are trying to prove, correspond to the sequence of unknowns $(e,j,g,h,t,v,w,z)$ as illustrated in Figure~\ref{secondFigure}\,.
\bigskip

We use the coordinates of the points on lines $\ell_1$ and $\ell_2$ to find equations of lines $A_1B_2, A_2B_1, A_1C_2, A_2C_1, B_1C_2$, and $B_2C_1$. This created the first difficulties, as the coefficients in the equations of these lines depended heavily on whether the corresponding parameters or unknowns were assumed to be in the basefield of the Hall system or not. Therefore we had to keep track of several possible equations for each of the six lines together with the assumptions on the nature of parameters and unknowns that were involved in finding each particular equation.

Then, using our choice of particular equations for each of the six lines, we found the coordinates of points $A_3, B_3$ and $C_3$ (if the corresponding pairs of lines intersected) or we showed that lines $A_1B_2, A_2B_1$, and $A_1C_2, A_2C_1$, and $B_1C_2, B_2C_1$ occur in parallel pairs.
While doing this, we again had to make assumptions on the coefficients of equations of lines being in or out of the basefield of the Hall system. When the points $A_3, B_3$ and $C_3$ exist, to check that there is indeed a Pappus configuration, we equated slopes of lines $A_3B_3$ and $A_3C_3$ or required that both slopes did not exist. This led to a system of equations with respect to the unknowns, and our goal was to prove that a solution exists. The number of various equations to solve grew extremely fast (it was in tens of thousands). As each line came with conditions on the parameters and unknowns that led to it, we tried to determined for each pair of lines whether these conditions were contradictory, and this allowed us to reduce the number of cases substantially (to about five thousand).

All this was done by using symbolic features of the software Mathematica \cite{mathematica} since the expressions for the coordinates of points $A_3$, $B_3$ and $C_3$ were unmanageable otherwise. As the symbolic power of Mathematica over finite fields is much weaker than it is over $\overline \Q$ (the default in Mathematica \cite{mathematica}), we attempted to solve the systems of equations symbolically assuming that all constants in our systems come from $\overline \Q$. In those cases, where Mathematica produced results over $\overline \Q$, it was generally easy to interpret them as outputs from symbolic computations over a finite field $\F=F_q$ for some prime power $q$. Some of the basic subroutines are listed in \cite{leshock}. Often Mathematica could not find results over $\overline \Q$, and we had to specialize some unknowns in a certain way to enable the program to find the remaining unknowns. Sometimes we could find this needed specialization of some unknowns and sometimes we could not. In the latter, we tried to argue that a solution existed.

In what follows, we present proof of the $2+0$ Theorem. For proof of Theorem \ref{theo31} that is largely similar to the one of the $2+0$ Theorem, we refer the reader to \cite{leshock}. Due to a limitation of space, we provide proof only for the cases for which the $2+0$ Theorem is not established in \cite{leshock}. Specifically, the proofs of the cases that use Propositions \ref{nbfnbfCases}\,, and \ref{nbfbfCases}\,(ii)\, are new. We present them below with the first case in greater detail.

\subsection{Case: $\boldsymbol{NBF/NBF \, (i)}$: Intersecting Lines}\label{subsectionNBFNBFIntersecting}

For this case, we consider any pair of intersecting lines from $NBF$.

\begin{proof}
We analytically construct a Pappus configuration to prove the $2+0$ Theorem. (This is also a proof that the $3+0$ Question is affirmative in this case.) Using Proposition \ref{nbfnbfCases}\, without loss of generality, we can assume that $\ell_1: \, \y = \x(\mu, \psi)$ and $\ell_2: \, \y=\x(0,1)$ for some $\mu, \psi \in \F$, $\psi \not=0$. We can determine the coordinates of points $A_1, B_1, C_1$ and $A_2, B_2, C_2$. Next, we build a type $2$ line $B_1C_2$ then use it to determine the equations of the remaining lines $A_1B_2, A_2B_1, A_1C_2, A_2C_1$ and $B_2C_1$. Each of the six type $2$ lines comes with constraints to ensure that the denominators and the second component of slope are nonzero. Then, we compute the the intersection point of lines $B_1C_2$ and $B_2C_1$, point $A_3$ and determine the conditions for its existence which we use to determine the conditions for the existence of and the coordinates of points $B_3$ and $C_3$. (Note: Our choice of which line or point to build first is arbitrary.)

In order for Mathematica to complete the computation that verifies that the slope of the line $A_3B_3$ equals the slope of the line $A_3C_3$ and therefore the Pappus line exists, we must carefully choose specializations of some unknowns that allow for a Pappus configuration to exist. This also serves to simplify the computation. Besides requiring that we have confirmation that the Pappus configuration exists, we must verify that each step in the computation is valid in the sense that we have not violated any laws of algebra or our own assumptions in any of the steps.

We begin with the previously created general formulas for the points $A_1, B_1, C_1$ and $A_2, B_2, C_2$. Using Proposition \ref{nbfnbfCases}\,, without loss of generality, we can fix the point $A_1$ by setting $\alpha=0$ and $\beta=1$. We provide a solution to this case by using specializations for the unknowns: $g=t=z=0$. The resulting equations of lines $\ell_1$ and $\ell_2$ and the coordinates of the points used in the Pappus configuration are listed below. Recall that $f(x)=x^2-rx-s$ is the defining polynomial of $\HH$.
\begin{talign*}
&\ell_1: \, (y_1,y_2)=(x_1,x_2)(\mu,\psi)&&\ell_2: \, (y_1,y_2)=(x_1,x_2)(0,1)\\
&A_1:\,((0, 1),(- f(\mu) / \psi, r-\mu)) &&A_2:\,((0,h),(s h, r h))\\
&B_1:\,((\gamma, \delta),(\gamma \mu - f(\mu) \delta / \psi,\delta(r-\mu)+\gamma \psi)) &&B_2:\,((0,v),(s v, r v))\\
&C_1:\,((e,j),(e \mu - f(\mu) j / \psi, j(r-\mu)+e \psi)) &&C_2:\,((w,0),(0, w))
\end{talign*}

Note that the denominator in some of the components of the points above is $\psi$ which is nonzero based on our assumption that $\ell_1$ is a type $2$ line. As a sample of the output in this case, the slope $\m$ and $\y$-intercept $\k$ of line $A_1B_2$ is provided below since it has the most compact formulas for any of the six lines $A_1B_2, A_2B_1, A_1C_2, A_2C_1, B_1C_2,$ and $B_2C_1$ that we can display in this case.

$m_1= \mu/(1-v)$ \\

$m_2 =(s (-1 + v)^2 - \mu (r (-1 + v) + \mu)) \psi /( (-1 + v) (\mu (-r + \mu) + s (-1 + v \psi)))$\\

$k_1 = v (s + (r - \mu) \mu - s \psi)/ ((-1 + v) \psi)$ \\

$k_2=v \mu / (1-v)$

For the denominators to be nonzero, we must not use $v=1$ or a $v$ that satisfies the quadratic equation in $v$, $s \psi \, v^2 +(-s + \mu(-r+\mu)-s \psi)\,v+s-\mu(-r+\mu) = 0$. To justify that this quadratic polynomial in $v$ is not identically zero for some sequence of parameters, consider the leading coefficient, $s \psi$. It is nonzero since $\psi \not=0$ and $-s$ is the constant term in the defining polynomial of $\HH$. Hence, in a large enough basefield, there exists a $v\not=1$ that is not a root of the quadratic polynomial (in $v$) listed above. Similar reasoning works for the denominators of the remaining five lines we build between the points on $\ell_1$ and $\ell_2$.

Since we chose to create type $2$ lines $A_1B_2$, $A_2B_1$, $A_1C_2$, $A_2C_1$, $B_1C_2$, and $B_2C_1$, we must also consider the constraints created to ensure that the second component of slope is nonzero for each line. We find that for line $A_1B_2$, the numerator of $m_2$ equals zero when $v=(2s+ r \mu \pm \sqrt{r^2+4s} \mu)/(2 s)$. We must avoid such values of $v$. There are many reasons why this does not invalidate our solution. First of all, if $\mu=0$, since $s \not=0$, we have that $v=1$ which we already determined that we would not use. If $\mu \not=0$, then these values of $v$ do not exist in the basefield due the fact that the defining polynomial $f(x)=x^2-rx-s$ is irreducible over $\F$. Similar reasoning ensures that the components of the other five lines listed above exist and the lines are of type $2$.

The intersection point of lines $A_1B_2$ and $A_2B_1$ is $C_3$. Its coordinates $(\x,\y)$ are provided below.

$x_1 =-(\gamma (v - h) h (\mu (-r + \mu) +
   s (-1 + v \psi)))/(-s (-v \delta + h)^2 + (v \delta \mu - h \mu -
    v \gamma \psi) (r (-v \delta + h) +
    v \delta \mu - h \mu - v \gamma \psi))$\\

$x_2=(s (\delta h^2 +
    v^2 \delta (\delta - h + \delta h) -
    v h (\delta + \delta^2 - h + \delta \
h)) - v^2 \delta^2 \mu^2 + v \delta h \mu^2 +
 v^2 \delta h \mu^2 + v \delta^2 h \mu^2 -
 v^2 \delta^2 h \mu^2 -
 v h^2 \mu^2 - \delta h^2 \mu^2 +
 v \delta h^2 \mu^2 +
 2 v^2 \gamma \delta \mu \psi -
 v^2 \gamma h \mu \psi -
 2 v \gamma \delta h \mu \psi +
 2 v^2 \gamma \delta h \mu \psi -
 v \gamma h^2 \mu \psi - v^2 \gamma^2 \psi^2 +
 v \gamma^2 h \psi^2 - v^2 \gamma^2 h \psi^2 +
 r (-\delta h +
    v (\delta - h + \delta h)) (v \delta \mu - \
h \mu - v \gamma \psi) )/(s (-v \delta + h)^2 + (v \delta \mu - h \mu -
    v \gamma \psi) (r v \delta - r h -
    v \delta \mu + h \mu + v \gamma \psi) )$\\

$y_1 =(-(v - h) (r - \mu) \mu (r \delta - \delta \mu + \
\gamma \psi) (v \delta \mu - h \mu -
    v \gamma \psi) -
 s^2 (v \delta - h) (-\delta h +
    v (\delta - h \psi + \delta h \psi)) +
 s (h^2 \mu (2 \delta \mu - \gamma \psi) +
    v^2 (\delta^2 \mu^2 (2 + h \psi) + \gamma \psi^2 \
(\gamma + h \mu + \gamma h \psi) - \delta \
\mu \psi (h \mu + 2 \gamma (1 + h \psi))) +
    r (-2 \delta h^2 \mu +
       v h (2 \delta^2 \mu - \gamma \delta \psi - \
h \mu \psi + \delta \mu (2 + h \psi)) +
       v^2 (-\gamma h \psi^2 - \delta^2 \mu (2 + \
h \psi) + \delta \psi (\gamma + h \mu + \
\gamma h \psi))) +
    v h (-2 \delta^2 \mu^2 - \delta \mu (-2 \gamma \
\psi + \mu (2 + h \psi)) + \psi (h \mu^2 - \
\gamma^2 \psi + \gamma (\mu + h \mu \psi)))) )/(\psi (-s (-v \delta + h)^2 + (v \delta \mu - h \
\mu - v \gamma \psi) (r (-v \delta + h) +
      v \delta \mu - h \mu - v \gamma \psi)) )$\\

$y_2 =(-r^2 (-\delta h +
    v (\delta - h + \delta h)) (v \delta \mu - \
h \mu -
    v \gamma \psi) - (v - h) (\mu (\delta \mu - \
\gamma \psi) (v \delta \mu - h \mu -
       v \gamma \psi) +
    s (-v \delta^2 \mu + \delta h \mu - \gamma \
h \psi + v \gamma h \psi)) +
 r (s (-\delta h^2 -
       v^2 \delta (\delta - h + \delta h) +
       v h (\delta + \delta^2 - h + \delta \
h)) + h^2 \mu (2 \delta \mu - \gamma \psi) +
    v^2 (\delta \mu - \gamma \psi) (\delta (2 + h) \
\mu - \gamma \psi - h (\mu + \gamma \psi)) +
    v h (-2 \delta^2 \mu^2 - \delta \mu ((2 + \
h) \mu -
          3 \gamma \psi) + \gamma \psi (\mu - \gamma \
\psi) + h \mu (\mu + \gamma \psi))) )/(-s (-v \delta + h)^2 + (v \delta \mu - h \mu -
    v \gamma \psi) (r (-v \delta + h) +
    v \delta \mu - h \mu - v \gamma \psi) )$

Note that the denominators of $x_1$ and $y_2$ are identical, and the denominator of $y_1$ is equal to $-\psi$ times the denominator of $x_2$. As stated above, the element $\psi \not=0$. The denominator of each component of point $C_3$ is a quadratic polynomial in $h$ with a nonzero leading coefficient, $-f(\mu)$ or $\psi f(\mu)$ (where $f$ is the defining polynomial of the Hall system which is irreducible over $\F$). Thus, in a large enough basefield, we can find values of $h$ that are not roots of those quadratic polynomials in $h$. In this way, we are certain that the components of point $C_3$ exist.

Finally, we create the lines $A_3B_3$ and $A_3C_3$ and compare the slopes to determine if the Pappus line exists. This solution yields a type $2$ Pappus line with a formula for the first component of slope that would fill $18$ pages of this document. It is too long to list here. The arguments for the existence of the components and that the lines are of type $2$ are similar to those used above. The subroutines to compute the components of the points and lines not listed here are provided on the website listed in citation \cite{leshock}\,.

\end{proof}

\subsection{Case: $\boldsymbol{NBF/NBF \, (ii)}$: Parallel Lines}\label{subsectionNBFNBFParallel}

For this case, we consider any pair of parallel lines from $NBF$.

\begin{proof}
We analytically construct a Pappus configuration to prove the $2+0$ Theorem. (This is also a proof that the $3+0$ Question is affirmative in this case.) We appeal to Proposition \ref{nbfnbfCases}\, to map any pair of parallel lines from $NBF$ to the pair $\ell_1: \, \y = \x(0,1) + \boldsymbol{\kappa}$ and $\ell_2: \, \y=\x(0,1)$ for some $\boldsymbol{\kappa} \in \HH$, and then we proceed as before to determine the coordinates of points $A_1, B_1, C_1$ and $A_2, B_2, C_2$ on lines $\ell_1$ and $\ell_2$, respectively. We build six type $2$ lines, $A_1B_2, A_2B_1, A_1C_2, A_2C_1, B_1C_2,$ and $B_2C_1$, then compute the coordinates of points $A_3$, $B_3$, and $C_3$.

Once again, we begin with general formulas for the points $A_1, B_1, C_1$ and $A_2, B_2, C_2$ and apply Proposition \ref{nbfnbfCases}\,, without loss of generality, to fix the point $A_1$ by setting $\alpha=0$ and $\beta=1$. We provide two solutions to this case by using specializations for the unknowns: (i) $g=h=t=0$ and (ii) $g=t=0, h=1$. The resulting equations of lines $\ell_1$ and $\ell_2$ and the coordinates of the points used in the Pappus configuration are listed below.
\begin{talign*}
(i)\,&\ell_1: \, (y_1,y_2)=(x_1,x_2)(0,1)+(\kappa_1, \kappa_2)&&\ell_2: \, (y_1,y_2)=(x_1,x_2)(0,1)\\
&A_1:\,((0, 1),(\kappa_1+s,\kappa_2+r)) &&A_2:\,((0,0),(0,0))\\
&B_1:\,((\gamma, \delta),(\kappa_1+s \delta,\kappa_2+\gamma+r \delta)) &&B_2:\,((0,v),(s v, r v))\\
&C_1:\,((e,j),(\kappa_1+s j,\kappa_2+e+rj)) &&C_2:\,((w,z),(sz, w+rz))
\end{talign*}

\begin{talign*}
(ii)\,&\ell_1: \, (y_1,y_2)=(x_1,x_2)(0,1)+(\kappa_1, \kappa_2)&&\ell_2: \, (y_1,y_2)=(x_1,x_2)(0,1)\\
&A_1:\,((0, 1),(\kappa_1+s,\kappa_2+r)) &&A_2:\,((0,1),(s,r))\\
&B_1:\,((\gamma, \delta),(\kappa_1+s \delta,\kappa_2+\gamma+r \delta)) &&B_2:\,((0,v),(s v, r v))\\
&C_1:\,((e,j),(\kappa_1+s j,\kappa_2+e+rj)) &&C_2:\,((w,z),(sz, w+rz))
\end{talign*}

The computations required to construct and to verify that the Pappus configuration exists proceed as before. We could not prove that the slopes of lines $A_3B_3$ and $A_3C_3$ were equal without specializing $h$. To show that we can specialize $h$ to only take on the values of $0$ or $1$ without loss of generality in the parameters, consider the lines $A_2B_1$ and $A_2C_1$ that both contain the point $A_2$. Since $A_2$ has component $x_2=h=0$ in solution (i) and $x_2=h=1$ in solution (ii), when we determine the constraints required so that the denominators and the second component of the slope of line $A_2B_1$ are nonzero, we find that if $\gamma=\kappa_2$, $\delta$ is arbitrary. If $\gamma=-\kappa_2$, then in solution (i), $\delta \not=0$ and in solution (ii) $\delta \not=1$. Since $\boldsymbol{\kappa}, \gamma$, and $\delta$ are given and the constraints are mutually exclusive, $\delta$ is arbitrary. Similarly, for line $A_1C_2$, we find that if $e =\kappa_2$, $j$ is arbitrary. If $e=-\kappa_2$, then in solution (i), $j \not=0$ and in solution (ii) $j \not=1$. Hence, $j$ is arbitrary. These are the only lines with constraints on parameters. All other possible constraints are resolved as they were in the previous case in Subsection \ref{subsectionNBFNBFIntersecting}. The verification that the polynomials in the denominators of the components of the lines containing point $A_2$ requires substituting in the constraints on the parameter sequence from one specialization to show that they are not constraints on the parameter sequence in the other specialization. The subroutines to compute the components of the points and lines not listed here are provided on the website listed in citation \cite{leshock}\,.
 \end{proof}

\subsection{Case: $\boldsymbol{BF/NBF}$}\label{subsectionBFNBF}

For this case, we consider any pair of lines with the first from $BF$ and the second from $NBF$. Clearly, these lines are intersecting.

\begin{proof}
We analytically construct a Pappus configuration to prove the $2+0$ Theorem. We appeal to Proposition \ref{nbfbfCases}\, to map any pair of lines with the first from $BF$ and the second from $NBF$ to the pair $\ell_1: \, \x = \0$ and $\ell_2: \, \y=\x(0,1)$, and then we proceed as before to determine the coordinates of points $A_1, B_1, C_1$ and $A_2, B_2, C_2$ on lines $\ell_1$ and $\ell_2$, respectively. We build six type $2$ lines, $A_1B_2, A_2B_1, A_1C_2, A_2C_1, B_1C_2,$ and $B_2C_1$, then compute the coordinates of points $A_3$, $B_3$, and $C_3$ or determine that lines $A_1B_2, A_2B_1$, and $A_1C_2, A_2C_1$, and $B_1C_2, B_2C_1$ occur in parallel pairs.

Once again, we begin with general formulas for the points $A_1, B_1, C_1$ and $A_2, B_2, C_2$ and apply Proposition \ref{nbfbfCases}\,, without loss of generality, to fix the point $A_1$ by setting $\alpha=0$ and $\beta=1$. We provide two solutions to this case by using specializations for the unknowns: If $\gamma=0$, then (i) $e=t=w=0$, and if $\gamma \not =0$, then (ii) $j=(\gamma - e+\delta e)/\gamma, \, t=~(\gamma^2+r \gamma g-s \delta g-s \gamma h)/(\gamma^2 +r \gamma \delta-s \delta^2), \, v=(\gamma \delta-\gamma g -s \delta h)/(\gamma^2+ r \gamma \delta-s \delta^2)$, $w=~(\gamma (s (-1 + \delta) e g+  \gamma \
(-e^2 + s g- r e g+
      s e h)))/(s ( \gamma + (-1 + \delta) e)^2 -  \gamma e (\
 \gamma e + r ( \gamma + (-1 + \delta) e))), z=(\gamma ( \gamma e (-1 + g) +
   s  \gamma h - (-1 + \delta) e (e -
      s h)))/(s ( \gamma + (-1 + \delta) e)^2 -  \gamma e (\
 \gamma e + r ( \gamma + (-1 + \delta) e)))$. The resulting equations of lines $\ell_1$ and $\ell_2$ and the coordinates of the points used in the Pappus configuration are listed below.
\begin{talign*}
(i)\,&\ell_1: \, (x_1,x_2)=(0,0)&&\ell_2: \, (y_1,y_2)=(x_1,x_2)(0,1)\\
&A_1:\,((0, 0),(0,1)) &&A_2:\,((g,h),(sh,g+rh))\\
&B_1:\,((0,0),(0, \delta)) &&B_2:\,((0,v),(s v, r v))\\
&C_1:\,((0,0),(0,j)) &&C_2:\,((0,z),(sz, rz))
\end{talign*}

\begin{talign*}
(ii)\,&\ell_1: \, (x_1,x_2)=(0,0)&&\ell_2: \, (y_1,y_2)=(x_1,x_2)(0,1)\\
&A_1:\,((0, 0),(0,1)
&&A_2:\,((g,h),(sh,g+rh))\\
&B_1:\,((0,0),(\gamma, \delta))
&&B_2:\, (\x_b,\y_b)\\
&C_1:\,\left(\left( 0,0  \right),\left( e, (\gamma-e+\delta e)/\gamma \right)\right)
&&C_2:\, (\x_c,\y_c)
\end{talign*}
where

$x_{b1}=(\gamma^2+r\gamma g-s \delta g-s \gamma h)/(\gamma^2+r \gamma \delta-s \delta^2)$

$x_{b2}=(\gamma \delta -\gamma g-s \delta h)/(\gamma^2+r \gamma \delta-s \delta^2)$

$y_{b1}=(s(\gamma(\delta-g)-s \delta h))/(\gamma^2+r \gamma \delta -s \delta^2)$

$y_{b2}=(\gamma^2 -s \delta(g+rh)+\gamma (r \delta-s h))/(\gamma^2 + r \gamma \delta-s \delta^2)$,

and

$x_{c1}=(\gamma (s (-1 + \delta) e g + \gamma \
(-e^2 + s g - r e g +
      s e h)))/(s (\gamma + (-1 + \delta) e)^2 - \gamma e (\
\gamma e + r (\gamma + (-1 + \delta) e)))$

$x_{c2}= (\gamma (\gamma e (-1 + g) +
   s \gamma h - (-1 + \delta) e (e -
      s h)))/(s (\gamma + (-1 + \delta) e)^2 - \gamma e (\
\gamma e + r (\gamma + (-1 + \delta) e)))$

$y_{c1}=(s \gamma (\gamma e (-1 + g) +
   s \gamma h - (-1 + \delta) e (e -
      s h)))/(s (\gamma + (-1 + \delta) e)^2 - \gamma e (\
\gamma e + r (\gamma + (-1 + \delta) e)))$

$y_{c2}= (\gamma ((-e^2 - e g r + g s + e h s) \gamma +
   e g s (-1 + \delta))+r \gamma (e (-1 + g) \gamma + h s \gamma -
   e (e - h s) (-1 + \delta)))/((s (\gamma + (-1 + \delta) e)^2 - \gamma e (\
\gamma e + r (\gamma + (-1 + \delta) e))))$

Both solutions have the six type $2$ lines, $A_1B_2, A_2B_1, A_1C_2, A_2C_1, B_1C_2,$ and $B_2C_1$. In solution (i), the Pappus line is also a type $2$ line. In solution (ii), the lines $A_1B_2, A_2B_1$, and $A_1C_2, A_2C_1$, and $B_1C_2, B_2C_1$ occur in parallel pairs. Let us list the reasons that these configurations together are sufficient to prove the $2+0$ Theorem in this case.

Solution (i) is used when $\gamma=0$ and the free unknowns on line $\ell_2$ are $g,h,v,z$. To verify that all denominators are nonzero and all type $2$ lines have slopes not in the basefield (nonzero second component) requires that the free unknowns can take on values that are not roots of certain polynomials. Using reasoning similar to that in the first case in Subsection \ref{subsectionNBFNBFIntersecting} to verify we have nonzero polynomials in a particular unknown, we find that in a large enough basefield, we are assured that there exist values of the unknowns that produce nonzero denominators and type $2$ lines with nonzero second component.

Solution (ii) is used when $\gamma \not=0$ and the free unknowns on line $\ell_2$ are $g,h$. To verify that all denominators are nonzero and all type $2$ lines have slopes not in the basefield (nonzero second component) requires that we ignore degenerate cases (when two of the points $A_1, B_1, C_1$ or $A_2, B_2, C_2$ coincide with each other or the point of intersection of lines $\ell_1, \ell_2$) which produce trivial Pappus configurations, that we recognize that the defining polynomial is irreducible over the basefield, and finally that we recognize that the free unknowns can take on values that are not roots of certain polynomials with at least one nonzero coefficient, which in a large enough basefield is assured. Note that in this case, the points $A_3, B_3$, and $C_3$ do not exist since the lines $A_1B_2, A_2B_1$, and $A_1C_2, A_2C_1$, and $B_1C_2, B_2C_1$ occur in parallel pairs.

The subroutines to compute the components of the points and lines not listed here are provided on the website listed in citation \cite{leshock}\,.
\end{proof}

\bigskip

\section{Concluding remarks}\label{sectionConcluding}

\rightline{{\it A fool can ask more questions in one hour than a wise man can answer in seven years.}}
\rightline{{\it -- European proverb.}}

\bigskip
As our verification using computer supports the affirmative answer to the $3+1$ Question, it is desirable to answer it (non-asymptotically) for Hall planes. Another line of research can be establishing the strongest form of Pappus's Theorem for other classes of finite nonclassical planes. We would like to repeat the question that motivated our research. We could not find any reference to it in the literature.

\begin{ques}Is it true that every finite projective plane contains at least one Pappus configuration?
\end{ques}

Let us call a pair of lines in a finite projective plane {\it Pappian} if any choice of three points on one line and any choice of three points on the other line yields a Pappus configuration. It was shown by Pickert \cite{Pickert}, and later by Burn \cite{Burn}, that any finite projective plane with a Pappian pair of lines is classical. We would like to ask the following question.

\begin{ques}\label{ques2} What is the smallest number of Pappus configurations on a pair of lines in a projective plane of order $n$ that implies the plane is classical?
\end{ques}

As was mentioned in the introduction, we know from \cite{ostromPaper} that every finite projective plane contains Desargues configuration.
Though this paper primarily concerns Pappus configuration, we think the following question, similar to Question \ref{ques2}\,, is of interest.
We call a triple of concurrent lines {\it Desarguesian} if any two triangles each having one vertex on each of these lines yields a Desargues configuration.

\begin{ques} (i) If a finite projective plane has a Desarguesian triple of lines, is it necessarily a classical plane?\\
\indent \indent \indent \indent (ii) What is the smallest number of Desargues configurations on a triple of concurrent lines in a projective plane of order $n$ that implies the plane is classical?
\end{ques}

Of course similar questions can be asked (and were asked) about some other configurations.
In particular, the one that attracted the attention of many researchers,
is the existence of a Fano configuration in any finite Non-Desarguesian
projective plane (see ~\cite{taitPaper}).

 At the end we want to mention an ``inverse'' question asked by Welsh \cite{Welsh} (see comments in \cite{moorhouseWillifordPaper}), and independently by Erd\H{o}s \cite{Erdos79} that in geometric terms can be stated as follows.

\begin{ques} Is every finite partial linear space (a configuration) embedded in a finite projective plane?
\end{ques}
See \cite{moorhouseWillifordPaper} for more details. As far as we know, researchers are divided in their opinions whether the answer to Erd\H{o}s' question is positive or negative. In graph theoretic terms, the question is equivalent to:
\begin{ques} Is every bipartite graph without 4-cycles a subgraph of the point-line incidence graph of a finite projective plane?
\end{ques}

\bigskip

\section{Acknowledgements}\label{sectionAcknowledgements}This work was partially supported by the Simons Foundation Award ID: $426092$ and the National Science Foundation Grant: $1855723$.

The authors are thankful to Eric Moorhouse for sharing his knowledge on Hall planes and, in particular, for correcting an error in the original description of the action of the collineation group on pairs of lines of the Hall plane. We are also thankful to Stefaan DeWinter, Bill Kantor, and Jason Williford
for useful discussions on the topics of this paper. Finally, we are thankful to the anonymous referees for their useful comments, and, in particular, for suggesting the references \cite{Burn}, \cite{Pickert}, and \cite{Welsh}.

\bigskip

\section{References}\label{sectionReferences}
\begingroup
\renewcommand{\section}[2]{}%

\endgroup

\end{document}